\documentclass{amsart}
\usepackage{amsmath,amsthm}
\usepackage{amsfonts,amssymb}

\usepackage{enumerate}

\usepackage{hyperref}

\usepackage{xcolor}
\usepackage[numeric]{amsrefs}

\newcommand{\assign}{:=}
\newcommand{\cdummy}{\cdot}
\newcommand{\nin}{\not\in}
\newcommand{\nobracket}{}

\newcommand{\tmop}[1]{\ensuremath{\operatorname{#1}}}
\newcommand{\tmscript}[1]{\text{\scriptsize{$#1$}}}

\newcommand{\tmtextrm}[1]{{\rmfamily{#1}}}

\theoremstyle{plain} 

\newtheorem{thm}{Theorem}[section]

\newtheorem{lem}[thm]{Lemma}

\newtheorem{exam}[thm]{Example}
\newtheorem{defn}[thm]{Definition}
\newtheorem*{thm*}{\bf Theorem}

\theoremstyle{remark}

\numberwithin{equation}{section}

\def\f{\frac}
\def\Bl{\Bigl}
\def\Br{\Bigr}
\def\dd{{\textnormal{d}}}
\def\ee{{\textnormal{e}}}
\def\i{{\textnormal{i}}}

\def\diam{\text{diam}}

\def\vi{\varphi}

\def\ta{\theta}
\def\al{{\alpha}}
\def\be{{\beta}}
\def\da{{\delta}}
\def\sa{{\sigma}}

\def\ga{{\gamma}}

\def\o{{\omega}}
\def\s{{\sigma}}
\def\va{\varepsilon}

\def\R{{\mathbb{R}}}

\def\ld{\lambda}

\def\EE{{\mathbb E}}

\def\NN{{\mathbb N}}

\def\RR{{\mathbb R}}
\def\SS{{\mathbb S}}

\def\dim{\operatorname{dim}}

\def\Re{\operatorname{Re}}
\def\Im{\operatorname{Im}}
\def\sph{\mathbb{S}^{d}}

\def\dsum{\displaystyle\sum}
\newcommand{\wt}{\widetilde}

\DeclareMathOperator{\pr}{Prob}

\DeclareMathOperator{\dist}{dist}
\def\Ld{\Lambda}
\def\Og{\Omega}

\input amssym.def

\begin{document}

\title[]{Discretization on high-dimensional domains}

  \author{ Martin D. Buhmann}

\address{Mathematics Institute,
	Justus-Liebig University\\
	Arndtstrasse 2\\
	D-35392 Giessen, Germany}
\email{buhmann@math.uni-giessen.de}

\author{Feng Dai}
\address{Department of Mathematical and Statistical Sciences\\
	University of Alberta\\ Edmonton, Alberta T6G 2G1, Canada.}
\email{fdai@ualberta.ca}

\author{Yeli Niu} 
\address{Department of Mathematical and Statistical Sciences\\
	University of Alberta\\ Edmonton, Alberta T6G 2G1, Canada.}
\email{yeli1@ualberta.ca}

\thanks{The second and the third  authors   were   supported
	by  NSERC  Canada under the
	grant  RGPIN 04702 Dai.
}


\begin{abstract} Let  $\mu$ be  a Borel probability measure on  a compact  path-connected  metric space $(X, \rho)$ for which there exist constants $c,\be>1$ such that  $\mu(B) \ge c r^{\be}$ 
	for every open ball $B\subset X$ of radius $r>0$.  For a class of Lipschitz functions $\Phi:[0,\infty)\to\RR$  that piecewisely lie  in a finite-dimensional  subspace of continuous functions, we prove under certain mild conditions on the metric $\rho$ and the  measure $\mu$  that  for each positive integer $N\ge 2$, and each  $g\in L^\infty(X, d\mu)$ with $\|g\|_\infty=1$,   there exist points  $y_1, \ldots, y_{ N}\in X$ and  real 
  numbers $\lambda_1, \ldots, \lambda_{ N}$  such that for any $x\in X$, 
  \begin{align*}
  &   \left| \int_X \Phi (\rho (x, y)) g(y) \,\dd \mu (y) - \sum_{j =
    1}^{ N} \lambda_j \Phi (\rho (x, y_j)) \right|   \leqslant C N^{- \frac{1}{2} - \frac{3}{2\be}}  \sqrt{\log N},
    \end{align*}
    where the constant  $C>0$ is  independent of $N$ and $g$. In the case when $X$ is the unit sphere $\sph$ of $\RR^{d+1}$ with the ususal geodesic distance, we also prove that the constant $C$  here is independent of the dimension $d$.  Our estimates are better than   those  obtained  from the standard  Monte Carlo methods,  which typically   yield  a weaker upper bound  $N^{-\f12}\sqrt{\log N}$. 

\end{abstract}

\maketitle

\section{Introduction}
\noindent
The theme of this article is the discretization in high-dimensional spaces and, using these discretizations, finding bounds for errors of numerical quadrature formulae. We mention the idea of the numerical
approximation of integrals by finite sums here at the beginning, because it is
a subject of interest generally in numerical analysis how to approximate the
integral of a function by a finite sum. The most basic approach to this is Gau\ss quadrature,
and the starting point of this in turn is using univariate numerical
integration employing zeros of orthogonal polynomials as knots \cite{Gautschi}.
The purpose of taking these zeros is to taylor the quadrature formula to provide
the optimal order of accuracy in approximating the integral by a finite sum.

The concept of Gau\ss quadrature (usually called cubature in higher
dimensions) can be generalized in many respects, see for instance \cite{FengDaiWang} for work on multivariate quadrature, and our goal in this
paper is to take a very general approach. To begin with, we shall work
in many (arbitrarily high) dimensions, and both in the literature and here,
multivariate spheres are of course our prime examples \cite{FengDaiYuanXu}, \cite{Petrova}, \cite{YuanXu}.

Secondly, we shall admit a
general metric space as the set on which our integrands are defined or
over which the integral shall be taken. We shall, third, find {\it
  dimension-independent\/} upper bounds on the error of cubature that
are uniform in $x$, where the integrals and sums take
the forms
\begin{align*}
  &   \int_X \Phi (\rho (x, y)) g(y) \,\dd \mu (y) 
\end{align*}
and\begin{align*}
  &   \sum_{j =
    1}^{ N} \lambda_j \Phi (\rho (x, y_j)),
    \end{align*}
respectively. Since these expressions depend on $x\in X$, where
$(X,\rho)$ is a compact metric space, they can be considered as a {\it
  discretization\/} of probability measures and $X$. It is attractive that the
upper bounds on the error are dimension-independent, because it allows us to
use these methods in high-dimensions without possibly large constants depending on dimensions marring our results.

With a constant $\beta$ depending on the Borel measure $\mu$, our goal
is to derive the estimate
 \begin{align*}
  &   \left| \int_X \Phi (\rho (x, y)) g(y) \,\dd \mu (y) - \sum_{j =
    1}^{ N} \lambda_j \Phi (\rho (x, y_j)) \right|   \leqslant C\|g\|_\infty  N^{- \frac{1}{2} - \frac{3}{2\be}}  \sqrt{\log N},
    \end{align*}
    where $C>0$ is a constant depending on $\Phi$ and certain
    properties of the measure $\mu$. The degrees of freedom to obtain
    the order of the estimate in $1/N$ on the right-hand side come
    from our judicious placement of the $y_j$s and the coefficients $\lambda_j$.  It is worthwile to point out here that such an estimate  is better than most  typical estimates that can be deduced from  the  Monte Carlo methods and 
    standard  probability techniques (based on various large deviation inequalities), which normally   yield  the  weaker upper bound  $N^{-\f12}\sqrt{\log N}$.

    Much of our work depends on the concepts of compactness, weak
    $(*)$ topologies and
    integrability with respect to a measure $\mu$, and therefore -- and
    for the purpose of fixing notation -- we shall review some of these
    points in the next section. 
    
    In Section~3, we  prove a preliminary
    result on the discretization of probability measures, which will play a vital role in this paper. To be more precise, let $Q$ be  a compact Hausdorff space equipped with     
     a Borel probability measure $\mu$, and let $X_m$ be an $m$-dimensional subspace of $C(Q)$.  Using the method of Bourgain and Lindenstrauss \cite{BL},   
      we   prove that for every $f\in C(Q)$,  the integral $\int_Q f(x)\, \dd\mu(x)$  can be discretized via  weighted sums 
    $$\sum_{j=1}^{m+2} \ld_j f(y_j),\  \ \ld_j\ge 0,\  \ y_j\in Q,$$
    where the weights $\ld_j\ge 0$ and the points $y_j\in Q$ are selected randomly according to a probability distribution in such a way that $\sum_{j=1}^{m+2} \ld_j=1$ and 
    $$\sum_{j=1}^{m+2} \ld_j f(y_j)=\int_Q f(x)\, \dd\mu(x),\  \ \forall f\in X_m.$$

    Section~4, then, considers regular partitions of compact metric
    spaces. Our main  result in this section, Theorem \ref{thm10},  states that  for a non-atomic  Borel probability measure $\mu$ on a compact path-connected metric space  $(\Omega, \rho)$ with diameter $\pi$,  there exists a partition $\{R_1, \ldots,
    R_N \}$ of $\Omega$ such that
     for each $1 \leq j \leq N$, $\mu (R_j) = \frac{1}{N}$ and
    	\text{\rm diam}$(R_j) \leq 4 \delta$, where $\da>0$ is a constant satisfying that 
    	 \begin{equation*}
    	\inf_{x \in \Omega} \mu \Bigl(B_{\delta / 2}^{} (x)\Bigr) \geqslant \frac{1}{N}.
    	\end{equation*}
       The crucial point here  lies in the fact that the constant $4$ in the estimate of $\diam(R_j)$  is absolute.

     Section~5 provides  one of the main results in
    Theorem~5.2.

    If the discretizations are to take place on {\it finite
      dimensional\/} compact domains, see, e.g., \cite{BD}, we have Theorem~6.2 as a
    suitable result.

    Sections~7 and 8 give some examples of interest, the example of the unit
    sphere probably giving the more important case, and Section~8
    suggesting some generalisations of the approach which is using piecewise
    polynomials $\Phi$ in our expressions above, to piecewise
    exponential functions instead.

\section{Preliminaries}

\noindent In this section, we list several basic results from functional analysis and
probability that will be needed in later sections. Most of the materials in
this section can be found in the book  {\cite{Wa}}.

\begin{thm}
  Let $X$ be a real linear topological space with dual space $X^{\ast}$. Then
  the following statements hold:
  \begin{enumerate}[\rm (i)]


    \item Let $A$ and $B$ be two nonempty disjoint convex sets in $X$.  If
    $A$ is open, then there exists $\Lambda \in X^{\ast}$ such that
    \[ \text{} \hspace{0.17em} \Lambda x < \inf_{y \in B}  \text{\tmtextrm{ }}
       \Lambda y,\quad \forall\; x \in A. \]
    If $A$ is compact, $B$ is closed and $X$ is locally convex, then there
    exist $\Lambda \in X^{\ast}$ such that
    \[ \sup_{x \in A}  \text{\tmtextrm{ }} \hspace{0.17em} \Lambda x < \inf_{y
       \in B}  \text{\tmtextrm{}} \hspace{0.17em} \Lambda y. \]

\item If   $X$ is an F-space ( i.e., a complete vector space
with metric that is translation invariant whose multiplications
and additions are continuous), then for every compact subset $K\subset X$,       the closure
of the  convex hull of $K$   is compact in $X$.
\end{enumerate}
  
\end{thm}

Next, we recall some basic facts on weak and weak*-topologies. A topology
$\tau_1$ on a nonempty set  $X$ is said to be weaker than another topology
$\tau_2$ on $X$ if $\tau_1 \subset \tau_2$.

\begin{thm}
  Let $X$ be a real  vector space, and   $X'$   a vector space of
  linear functionals on $X$ which separates points in $X$ (i.e.,  given  any two distinct
  points $x_1, x_2 \in X$ there exists $\Lambda \in X'$ such that $\Lambda x_1
  \neq \Lambda x_2$).  If  $\tau$  denotes  the weakest topology on $X$ with
  respect to which every  element in $X'$ is a continuous linear functional
  on $X$, then $(X, \tau)$ is a locally convex space whose dual
  is $X'$.
\end{thm}

Let  $X$ be  a real, locally convex  linear  topological  space with topology $\tau$ and the dual space   $X^{\ast}$.  Let   $\tau_w$ denote the weak topology of $X$, i.e.,   the weakest
topology of $X$ with respect to  which every linear functional  in $X^{\ast}$ is
continuous.  Then $\tau_w\subset \tau$, and  $X_w=(X,\tau_w)$  is  a locally convex space whose dual is also
  $X^{\ast}$.
We denote by $\tau_{w^{\ast}}$ the weak${}^\ast$ -topology of
  $X^{\ast}$; that is,  $\tau_{w^{\ast}}$ is the weakest
  topology of $X^{\ast} $ with respect to which for  every $x\in X$, the linear functional $f\in X^\ast\to f(x)$  is continuous. Then $(X^{\ast}, \tau_{w^{\ast}})$ is a locally convex linear topological
  space whose dual is $X$. If  $X$ is separable, then every weak*-compact set
  $K$ in $X^{\ast}$  is metrizable in the weak*-topology.

\begin{thm} [Banach-Alaoglu theorem] For every neighborhood $V$ of $0$ in $X$, its polar 
	\[ K \assign \{\Lambda \in X^{\ast} : | \Lambda x| \leq 1, \forall x \in
	V\} \]
	is weak* -compact in $X^\ast$.  If, in addition, $X$ is separable, then $K$  is sequentially compact in the weak* -topology.
\end{thm}

Third, we review some basic results on vector-valued integration. We
start with the following definition:

\begin{defn}
  Let $X$ be  a  real locally convex  topological vector space, and let  $(Q, \mu) $ be a measure space. A vector-valued  function $f : Q
  \to X$ is said to be integrable with respect to $\mu$ if
  \[ \Lambda (f (\cdot)) = \langle \Lambda, f (\cdummy) \rangle \in L^1 (Q,
     \mu), \qquad \forall \Lambda \in X^{\ast} \quad \]
  and there exists $y \in X$ such that
  \[ \langle \Lambda, y \rangle = \int_Q \langle \Lambda, f (x) \rangle \,\dd \mu
     (x), \qquad \forall \Lambda \in X^{\ast} . \]
  If such a vector $y \in X$  exists, it must be unique, and is denoted
  by  $\int_Q f (x) \,\dd  \mu (x)$.
\end{defn}

Recall that a positive Borel measure $\mu$ on a topological space  $Q$ is
regular if

\begin{align*}
  \mu (E) & = \sup \{\mu (K) : K \subset E \text{ is compact} \}\\
  & = \inf \{\mu (G) : E \subset G, G \text{ is open in } X\}
\end{align*}
for every Borel set $E \subset Q$.\quad Each \ Borel probability measure on a
locally compact Hausdorff space with a countable base for its topology, or on
a  compact metric space is regular. If $Q$ is a  compact Hausdorff space,
and $C_{} (Q)$ is the space of all continuous functions on $Q$ (with the
uniform norm), then the dual of $C (Q)$ is the space of all finite
regular Borel measures (i.e., Radon measures) on $Q$ (with the norm of total
variation).

\begin{thm}
  \label{thm-2-2-0}Suppose that
  \begin{enumerate}[\rm (i)]
    \item $X$ is a real, locally convex  topological vector space;

    \item $Q$ is a compact Hausdorff space;

    \item $f : Q \to X$ is continuous;

    \item $\overline{\tmop{conv} (f (Q))}$ is compact in $X$ (this is
    automatically true if $X$ is an F-space).
  \end{enumerate}
  Then given any Borel probability measure  $\mu$ on $Q$, the function $f : Q \longrightarrow X$ is integrable with respect to $\mu$  and
  moreover,
  \[ y = \int_Q f \hspace{0.17em} \,\dd \mu = \int_{f (Q)} z \hspace{0.17em}
  \,\dd
     \mu_f (z) \in \overline{\tmop{conv} (f (Q))}, \]
     where $\mu_f$ is a Borel probability measure on $f(Q)$ given by 
     $$ \mu_f (E)=\mu(f^{-1}(E)),\   \ E\subset f(Q).$$

  Conversely, if $y\in \overline{\tmop{conv} (f (Q))}$, then there exists a regular Borel probability measure $\mu_f$ on $f(Q)$ such that 
   \[ y = \int_{f (Q)} z \hspace{0.17em}
  \,\dd
  \mu_f (z).\]
\end{thm}


\begin{thm}
  Suppose that $Q$ is a compact Hausdorff space, $X$ is a Banach space, $f : Q
  \to X$ is continuous, and $\mu$ is a positive Borel measure on $Q$. Then
  \[ \left\| \int_Q f \,\dd  \mu \right\| \leq \int_Q \|f\|
     \,\dd  \mu . \]
\end{thm}

\section{A preliminary result}

\noindent Let $Q$ be a compact metric space equipped with a Borel probability measure
$\mu$. Let $M (Q)$ denote the space of all finite signed Borel measures on
$Q$. Then $M (Q)$ is a Banach space with respect to the norm
\[ \| \nu \| \assign | \nu | (Q) = \sup \left\{ \left| \int_Q f\,\dd \nu \right| :
   f \in C (Q), \|f\|_{C (Q)} \leq 1 \right\} . \]
Such a Banach space is the dual space of $C (Q)$. Note that $C (Q)$ is a
separable Banach space.  Let $M (Q)^{w^{\ast}}$ denote the space $M (Q)$
endowed with the weak* -topology $\tau_{w^{\ast}}$. Then $M (Q)^{w^{\ast}}$ is
a locally convex topological space with dual space $C (Q)$.

Next, let $X_m$ denote an $m$-dimensional linear subspace of $C (Q)$. Let
$\Sigma_0 \subset M (Q)$ denote the set of all probability measures $\rho \in
M (Q)$ of the form $$\rho = \sum_{j = 1}^{m + 2} \lambda_j (\rho) \delta_{y_j
(\rho)},$$ where $\lambda_j (\rho) \ge 0$, $y_j (\rho) \in Q$ for $j = 1, 2,
\ldots, m + 2$ and $\sum_{j = 1}^{m + 2} \lambda_j (\rho) = 1$.

Let $\Sigma
\subseteq \Sigma_0$ denote the set of all probability measures $\rho \in
\Sigma_0$ such that
\[ \int_Q f (x)  \,\dd  \mu (x) = \int_Q f (x)\,\dd \rho (x), \quad\forall
   f \in X_m . \]
\begin{thm}
  \label{thm-2-1}There exists a Borel probability measure $\nu$ on the space
  $M (Q)^{w^{\ast}}$ which is supported in the set $\Sigma \subset M (Q)$ and
  satisfies
  \[ \mu = \int_{\Sigma} \rho\,\dd  \nu (\rho), \]
  where the equality holds in the sense that for any $f \in C (Q)$,
  \[ \int_Q f (x) \,\dd \mu (x) = \int_{\Sigma} \sum_{j = 1}^{m + 2} \lambda_j
  (\rho) f (y_j (\rho))  \hspace{0.17em} \dd \nu (\rho)  \]
  and where $\mu$ is the probability measure we wish to discretise.
\end{thm}

\begin{lem}
  \label{lem-6-2}The set $\Sigma$ is $w^{\ast}$-compact in $M (Q)$.
\end{lem}

\begin{proof}
  Define
  \[ S \assign \left\{ \lambda = (\lambda_1, \ldots, \lambda_{m + 2}) \in
     \mathbb{R}^{m + 2} : \lambda_1, \ldots, \lambda_{m + 2} \ge 0, \sum_{j =
     1}^{m + 2} \lambda_j = 1 \right\} . \]
  Then $S \times Q^{m + 2}$ is a compact topological space with respect to the
  product topology. Next, consider the mapping $T : S \times Q^{m + 2} \to M
  (Q)^{w \ast}$ that takes $(\lambda, x) \in S \times Q^{m + 2}$ to the
  measure $\sum_{j = 1}^{m + 2} \lambda_j \delta_{x_j} \in M (Q)$. Note that
  for any $f \in C (Q)$, and any $(\lambda, x), (\alpha, y) \in S \times Q^{m
  + 2}$, we have

  \begin{align*}
    \left| \Big\langle \sum_{j = 1}^{m + 2} \lambda_j \delta_{x_j} - \sum_{j =
    1}^{m + 2} \alpha_j \delta_{y_j}, f \Big\rangle \right| & \leq \sum_{j = 1}^{m
    + 2} | \lambda_j f (x_j) - \alpha_j f (y_j) |\\
    & \to 0, \text{as } (\alpha, y) \to (\lambda, x) .
  \end{align*}

  This implies that the mapping $T$ is continuous, and hence $\Sigma_0 = T (S
  \times Q^{m + 2})$ is $w^{\ast}$-compact.

  Finally, for each $f \in C (Q)$, set $\mu_f \assign \int_Q f \hspace{0.17em}
  \,\dd \mu$. Then
  \[ \Sigma \assign \{\rho \in \Sigma_0 : \langle f, \rho \rangle = \mu_f,
     \forall f \in X_m \} . \]
  Since each $X_m \subset C (Q)$ and $C (Q)$ is the dual space of $M
  (Q)^{w^{\ast}}$, it follows that $\Sigma$ is a $w^{\ast}$-closed subset of
  the $w^{\ast}$-compact set $\Sigma_0$. Thus, $\Sigma$ is a weak*-compact
  subset of $M (Q)$.
\end{proof}

\begin{lem}
  \label{lem-6-3}The probability measure $\mu \in M (Q)$ is in the
  weak*-closure of the convex hull $K$ of $\Sigma \subset M (Q)^{w^{\ast}}$.
\end{lem}

\begin{proof}
  Assume to the contrary that \ $\mu \nin K = \overline{\mathrm{conv}
  \hspace{0.17em} \Sigma}^{^{^{w^{\ast}}}}$.\quad Then \ by the convex
  separation theorem, there exists $g \in C (Q)$ such that
  \begin{equation}
    \label{2-1} \int_Q g \,\dd  \mu > \sup_{\rho \in \Sigma}
    \int_Q g \,\dd  \rho .
  \end{equation}
  Let $X_{m + 1} = \tmop{span} \{ X_m, g \} .$ By Corollary 4.1 of \cite{DPTY}, there
  exist $x_1, x_2, \ldots, x_{m + 2} \in Q$\quad and $\lambda_1, \ldots,
  \lambda_{m + 2} \geqslant 0$ such that $\sum_{j=1}^{m+2} \ld_j=1$ and 
  \[ \int_Q f \,\dd \mu = \sum_{j = 1}^{m + 2} \lambda_j f (x_j), \quad \forall f
     \in X_{m + 1} . \]
  This implies that $\rho = \sum_{j = 1}^{m + 2} \lambda_j \delta_{x_j} \in
  \Sigma$ and $\int_Q g \,\dd \mu = \int_Q g d \rho$, which contradicts
        {\eqref{2-1}}.

\end{proof}

\subsection*{Proof of Theorem \ref{thm-2-1}}Let $X = C (Q)$. Then $M (Q) =
X^{\ast} .$ By Lemma \ref{lem-6-3}, \ $\mu$ lies in the $w^{\ast}$-closure  of
the convex hull of $\Sigma$; that is, $\mu \in K \assign \overline{\tmop{conv}
(\Sigma)}^{^{^{w^{\ast}}}} .$ By Lemma \ref{lem-6-2} , \ $\Sigma$ is compact
in the space \ $(X^{\ast}, w^{\ast})$. Thus, by Theorem \ref{thm-2-2-0}, it is
enough to show that $K$ is also compact in the space $(X^{\ast},
w^{\ast})$.\quad Note that
\[ \Sigma \subset \Sigma_0 \subset B_{X^{\ast}} \assign \{ \nu \in X^{\ast} :
   \| \nu \| \leqslant 1 \}, \]
which also implies that $\tmop{conv} (\Sigma) \subset B_{X^{\ast}} $. Since
$B_{X^{\ast}}$ is compact in the space $(X^{\ast}, w^{\ast})$, it follows that
$K \assign \overline{\tmop{conv} (\Sigma)}^{^{^{w^{\ast}}}}$is a closed subset
of $B_{X^{\ast}}$, which also implies that $K$ is compact in the space
$(X^{\ast}, w^{\ast})$. The theorem is proved.

\section{Regular partitions on compact metric space}
\noindent Let $(\Omega, \rho)$ be a compact metric space$.$ Open balls and closed balls  in $\Omega$ will
be denoted by $B_\zeta (x) \assign \{ y \in \Omega : \rho (x, y) < \zeta \}$, and
$B_\zeta [x] \assign \{ y \in \Omega : \rho (x, y) \leqslant \zeta \}$, respectively.
 A path
connecting two points $x, y \in \Omega$ is a continuous map $\gamma : [0, 1]
\rightarrow \Omega$ with $\gamma (0) = x$ and $\gamma (1) = y$. A metric space $(\Og, \rho)$ is called path-connected if every two distinct points in $\Og$ can be connected with a path. As is well known, every open connected subset of $\RR^n$ is path-connected. 
 Given a
set $A \subset \Omega$ and a point \ $x \in \Omega$, define
\[ \tmop{dist} (x, A) \assign \inf_{y \in A} \rho (x, y) . \]

\begin{thm}
  \label{thm10}Let $(\Omega, \rho)$ be a compact path-connected  metric space with
  diameter $\tmop{diam} (\Omega) \assign \max_{x, y \in \Omega} \rho (x, y) =
  \pi .$ Let  $\mu$ be  a non-atomic  Borel probability measure on
  $\Omega$, and  $N\ge 2$  a positive integer. Assume that the inequality
  \begin{equation}
    \inf_{x \in \Omega} \mu \Bigl(B_{\delta / 2}^{} (x)\Bigr) \geqslant \frac{1}{N}
    \label{3-3}
  \end{equation}
  holds for some $\delta > 0$. Then there exists a partition $\{R_1, \ldots,
  R_N \}$ of $\Omega$ such that
  \begin{enumerate}[\rm (i)]
    \item the $R_j$ are pairwise disjoint subsets of $\Omega$,

    \item for each $1 \leq j \leq N$, $\mu (R_j) = \frac{1}{N}$ and
    \text{\rm diam}$(R_j) \leq 4 \delta$. 
  \end{enumerate}
\end{thm}

    Theorem \ref{thm10} with constants depending on certain geometric  parameters of the underlying  space $(\Og, \rho,\mu)$ (e.g. {\it dimension, doubling constants}) is probably known in a more general setting. The crucial point here lies in the fact that the constant $4$ in the estimates of $\diam(R_j)$ is absolute.

\begin{lem}
  \label{lem-1-2} Let $(\Omega, \rho)$ be a compact path-connected  metric space with diameter $\pi$.
  Then for each $\delta \in (0, \pi)$, there exist a finite set $\Lambda =
  \{a_1, \ldots, a_M \} \subset \Omega$ with $M>1$ such that $\Omega = \bigcup_{j = 1}^M
  B_{\delta} (a_j)$ and
  \[ \tmop{dist} (a_j, \Lambda_{j - 1}) = \delta, \quad j = 2, 3, \ldots, M,
  \]
  where $\Lambda_k \assign \{ a_1, a_2, \ldots, a_k \}, k = 1, \ldots, M.$
\end{lem}

\begin{proof}
   Since the metric space $\Og$ is path-connected and has diameter $\pi\ge \da$, there exist two  points $a_1, a_2\in\Og$ such that $\rho(a_1, a_2)=\da$.
  Assume that
  $\Lambda_n = \{ a_1, \ldots, a_n \}$ is a finite subset of $\Omega$ such that
  \[ \tmop{dist} (a_j, \Lambda_{j - 1}) = \delta, j = 2, \ldots, n, \]
  where $\Lambda_j = \{ a_1, a_2, \ldots, a_j \}.$
   If $\Omega = \bigcup_{j = 1}^n B_{\delta} (a_j)$, then it is sufficient to use $M = n$. Now assume that, in contrast, $\Omega \neq \bigcup_{j = 1}^n
  B_{\delta} (a_j)$. Then  there exists a point  $y \in \Omega \setminus \Lambda_n$ such
  that $$\tmop{dist} (y, \Lambda_n) \geqslant \delta.$$
  Without loss of generality, we may assume that $\text{dist} (y, \Lambda_n)=\rho(y, a_1)$.
  Let $\ga:[0,1]\to \Omega$ be a path such that $\ga(0)=y$ and
  $\ga(1)=a_1$. Define  $f(t):=\dist(\ga(t),\Ld_n)$ for $t\in
  [0,1]$. Clearly, $f$ is a continuous function on $[0,1]$ with \[
  f(0) =\dist(y, \Ld_n)\ge \da\   \  \text{and}\   \ f(1)
  =\text{dist}(a_1, \Ld_n)=0.\]
  Thus, there exists a point $a_{n+1}=\ga(t_n)\in \Og$ for some $t_n\in [0, 1]$ such that $$\text{dist} (a_{n+1},\Ld_n)=f(t_n)=\da.$$
  We may continue this selection procedure with $\Lambda_{n + 1} = \{a_1,
  \ldots, a_{n + 1} \}$. Since $\Omega$ is compact, this procedure must
  terminate after a finite number of steps.
\end{proof}

\begin{proof}[Proof of Theorem \ref{thm10}]
  Let $$\{a_1, \ldots, a_M \}$$ be a finite subset of $\Omega$ as given in
  Lemma~ \ref{lem-1-2}.

  For $1 < j \leq M$, let $1 \leq k_j < j$ be an integer
  such that
  \[ \tmop{dist} (a_j, \Lambda_{j - 1}) = \rho (a_j, a_{k_j}) = \delta. \]
  For each $1 \leq j \leq M$, define
  \[ V_j \assign \left\{ x \in \Omega : \rho (x, a_j) = \dist (x,
     \Lambda)\  \  \text{and}\  \  \dist (x, \Lambda) < \min_{1 \leq i < j}
     \rho (x, a_i) \right\} . \]
  That is, $x \in V_j$ if and only if $j$ is the smallest positive integer  such that $\dist (x,
  \Lambda) = \rho (x, a_j).$ Clearly,  the sets $V_j$ are pairwise disjoint,
  \begin{equation}
    \label{1-3} B_{\frac{\delta}{2}} (a_j) \subset V_j \subset B_{\delta}
    [a_j],\   \  j = 1, 2, \ldots, M,
  \end{equation}
  and $\Omega = \bigcup_{j = 1}^M V_j$. Moreover, using {\eqref{3-3}}, we have
  \[ \mu (V_j) \geqslant \frac{1}{N}, \quad \forall 1 \leqslant j \leqslant
     M. \]

  Now we  construct the desired  partition of $\Omega$ as follows  via a finite number
  of steps.
   In the first step, we write $V_j^0 = V_j$ for $j = 1, \ldots, M$, and
    modify the cells $V_M$ and $V_{k_M}$ slightly so that
  $N \mu (V_M) $ is an integer. Let $E_M \subset V_M^0$ be such that $\mu
  (E_M) < \frac{1}{N}$ and $N \mu (V_M^0 \setminus E_M)$ is a positive
  integer. We then update the cells as follows:
  \[ V_j^1 \assign \left\{\begin{array}{ll}
       V_j^0, & \text{if $j \neq M$ and } j \neq k_M ,\\
       V_j^0 \setminus E_M, & \text{if } j = M,\\
       V_j^0 \cup E_M, & \text{if } j = k_M .
     \end{array}\right. \]
  Note that the sets $V_j^1$ are pairwise disjoint, $\Omega = \bigcup_{j =
  1}^M V_j^1$, $V_j^0 \subset V_j^1$ for $1 \leq j \leq M - 1$ and $V_M^1
  \subset V_M^0$.

  In the second step, we continue the process  with the collection of the first $M-1$ updated cells: $V_j^1$, $1\leq j\leq M-1$. More precisely,   we choose a subset  $E_{M - 1}$ of  $ V_{M - 1}^0$ such that $\mu (E_{M - 1}) < \frac{1}{N}$ and $N \mu (V_{M-1}^1
  \setminus E_{M - 1})$ is a positive integer, and then   update the cells as
  follows:
  \[ V_j^2 \assign \left\{\begin{array}{ll}
       V_j^1, & \text{if $j \neq M - 1$ and } j \neq k_{M - 1},\\
       V_j^1 \setminus E_{M - 1}, & \text{if } j = M - 1,\\
       V_j^1 \cup E_{M - 1}, & \text{if } j = k_{M - 1} .
     \end{array}\right. \]
      It is very important here that the set $E_{M-1}$ is selected  as a subset of $V_{M-1}^0$ (rather than  a general subset $V_{M-1}^1$)  because this  way of selection yields  a better  control of  the diameter of the updated cell   $V_{k_{M-1}}^1:=E_{M-1}\cup V_{k_{M-1}}^1$.

  In general, \ at the $\ell$-th step with $1 \leq \ell < M$, we modify the
  cells $V_{M - \ell + 1}^{\ell - 1}$ and $V_{k_{_{M - \ell + 1}}}^{\ell - 1}$ in
  a similar manner. Indeed, let $E_{M - \ell + 1} \subset V_{M - \ell + 1}^0
  \subset V_{M - \ell + 1}^{\ell - 1}$ be such that $\mu (E_{M - \ell + 1}) <
  \frac{1}{N}$\quad and\quad$N \mu (V_{M - \ell + 1}^{\ell - 1} \setminus E_{M -
  \ell + 1})$ is a positive integer.\quad We then \ define
  \[ V_j^{\ell} \assign \left\{\begin{array}{ll}
       V_j^{\ell - 1}, & \text{if $j \neq M - \ell + 1$ and } j \neq k_{M -
       \ell + 1},\\
       V_{M - \ell + 1}^{\ell - 1} \setminus E_{M - \ell + 1}, & \text{if } j
       = M - \ell + 1,\\
       V_{k_{M - \ell + 1}}^{\ell - 1} \cup E_{M - \ell + 1}, & \text{if } j =
       k_{M - \ell + 1} .
     \end{array}\right. \]
  Clearly, the sets $V_j^{\ell}$ are pairwise disjoint, $\Omega = \bigcup_{j =
  1}^M V_j^{\ell}$,
  \[ V_j^0 \subset V_j^{\ell - 1} \subset V_j^{\ell} \quad \tmop{for} \quad j
     = 1, 2, \ldots, M - \ell, \]
  and for $j = M - \ell + 1, \ldots, M$,
  \[ V_j^{\ell} \subset V_j^{\ell - 1}  \text{\ \ \ and \ $N \mu (V_j^{\ell})$ is a
     positive integer.} \]
  Furthermore, by the above construction, it is easily seen that \ for each $1 \leq j \leq
  M - \ell$,
  \[ V_j^{\ell} \subset \bigcup_{\tmscript{\begin{array}{c}
       M - \ell + 1 \leq k \leq M\\
       \rho (a_k, a_j) = \delta
     \end{array}}} (V_j^0 \cup V_k^0) , \]
 which, using {\eqref{1-3}}, implies that $V_j^{\ell} \subset B_{2
  \delta}  [a_j]$ and
  $ \text{diam} (V_j^{\ell}) \leq 4 \delta$ for all $1 \leq j \leq M$.

  The  above process will be terminated  after the  $ (M - 1)$-st step, where we obtain  pairwise disjoint subsets $V_j^{M - 1}$, $j = 1, 2, \ldots, M,$ of
  $\Omega$ with diameter $\leq 4 \delta$ such that $\Omega = \bigcup_{j = 1}^M
  V_j^{M - 1}$ and $N \mu (V_j^{M - 1})$ is a positive integer for $2 \leq j
  \leq M$. Since $\mu$ is a probability measure, we have
  \[ N = N \mu (\Omega) = \sum_{j = 1}^M N \mu (V_j^{M - 1}) . \]
  This implies \ that $N \mu (V_1^{M - 1})$ is a positive integer as
  well.
  Since $\mu$ is non-atomic,  for each $1 \leq j \leq M$, we may write
  $V_j^{M - 1}$ as a disjoint union
  \[ V_j^{M - 1} = \bigcup_{k = 1}^{\ell_j} S_{j, k} \]
  such that $\mu (S_{j, k}) = \frac{1}{N}$ and diam$(S_{j, k}) \leq 4 \delta$
  for $1 \leqslant k \leqslant \ell_j$. This leads to a partition of $\Omega$
  with the desired properties:
  \[ \Omega = \bigcup_{j = 1}^M \bigcup_{k = 1}^{\ell_j} S_{j, k} . \]
\end{proof}

\section{Discretization  on compact  metric spaces }

\noindent
Let $(X, \rho) $ be a compact   metric space with metric $\rho$ and
diameter $ \pi$. For $x\in X$ and $0\leq  \underline a< \underline
b\leq \pi$, set  
$$E(x; \underline  a, \underline  b):=\{y\in X:\  \ \underline a\leq
\rho(x,y)\leq  \underline b\}.$$ 
A partition of $X$ consists of finitely many pairwise disjoint subsets
of $X$ whose union is $X$. 
\begin{defn} 
	  Let
	$0 = t_0 < t_1 < \cdots < t_{\ell} = \pi$  be a partition of  the interval
	$[0, \pi]$, and let $r\in\NN$.  We say  $\Phi\in C[0,\pi]$
          belongs to  the class    $\mathcal{S}_r\equiv
          \mathcal{S}_{r}(t_1,\ldots, t_\ell) $ if  there exists an
          $r$-dimensional  linear subspace $V_r$ of $C(X)$ such that
          for any $x\in X$ and each $1\leq j\leq \ell$,   
	$$\Phi(\rho(x,\cdot))\Bl|_{E(x; t_{j-1}, t_j)}\in\Bl\{
          f\Bl|_{E(x; t_{j-1}, t_j)}:\  \ f\in V_r\Br\}.$$ 
	\end{defn}

Next, let  $\mu$ be  a Borel probability measure on $X$ satisfying the following condition for a parameter  $\be\ge 1$ and some  constant $c_1>1$:\\

\begin{enumerate}[\rm (a)]
	\item for each positive integer $N$,  there exists a partition $\{X_1, \ldots,
	X_N \}$ of $X$ such that    $\mu (X_j) = \frac{1}{N}$ and
	\text{diam}$(X_j) \leq  \delta_N:=c_1N^{-\f 1\be}$ for $1\leq j\leq N$.

\end{enumerate}

According to  Theorem \ref{thm10},  Condition (a)   holds automatically with $c_1= 20\pi$  if the metric space  $X$ is path-connected,  and  $\mu$ is a non-atomic Borel probability measure on $X$ satisfying that   for any $0<t\leq 1$, 
\begin{equation}\label{4-1}
\inf_{x\in X}  \mu (B_t (x)) \ge \Bl( \f 8 {c_1}\Br)^\be t^\be.
\end{equation}

In this section, we shall prove 

\begin{thm}
	\label{thm15}  	
	Let $\Phi\in C[0,\pi] $
	satisfy 	\begin{equation}
	\label{4-3-0} | \Phi (s) - \Phi (s') | \leq |s - s'|,\qquad   \forall s, s' \in [0,
	\pi],
	\end{equation}
	and  belong to  a   class   $\mathcal{S}_{r}(t_1,\ldots,
        t_\ell) $ for some compact metric space $(X,\rho)$, where
        $r\in \NN$  and 	$0 = t_0 < t_1 < \cdots < t_{\ell} =
        \pi$. 	Let  $\mu$ be  a Borel probability measure on  $X$
        satisfying the    condition (a)  and  
 the following condition:\vspace{2mm}
	
	\begin{enumerate}
		
	\item [{\rm (b)}]
	 for each $x \in X$  and $\da\in (0, \pi)$,
	\begin{equation}\label{4-1-0}
	\mu \Bl(  E(x; t_j-\da, t_j+\da)\Br) \leqslant c_2 \delta,\   \ 1\leq j<\ell,
	\end{equation}
	where $c_2>1$ is a constant independent of $\da$ and $x$. 
	\end{enumerate}\vspace{2mm}
	Then for each positive integer $N\ge 4$ 
	, there exist points  $y_1, \ldots, y_{_{(r + 2) N}} \in X$ and  nonnegative
	numbers $\lambda_1, \ldots, \lambda_{_{(r + 2) N}} \geqslant 0$ such that $
	\dsum_{j = 1}^{(r + 2) N} \lambda_j = 1$ and
	\begin{align*}
	&  \max_{x \in X} \left| \int_X \Phi (\rho (x, y)) \,\dd \mu (y) - \sum_{j =
		1}^{(r + 2) N} \lambda_j \Phi (\rho (x, y_j)) \right|   \leqslant c_3  N^{- \frac{1}{2} - \frac{3}{2\be }}  \sqrt{\log N},
	\end{align*}
	where $c_3: =8 c_1^{2} \sqrt{c_2\ell} \sqrt{\be}$.
\end{thm}



In the case when the metric space $X$ is  path-connected, we will prove

\begin{thm}\label{cor-4-2} Let $(X,\rho)$ be a compact path-connected metric space. 	Let $\Phi\in C[0,\pi] $ satisfy  \eqref{4-3-0} and belong to  a  class   $\mathcal{S}_{r}(t_1,\ldots, t_\ell) $ for some 
	  $r\in \NN$ and  	$0 = t_0 < t_1 < \cdots < t_{\ell} = \pi$. Let  $\mu$ be  a non-atomic Borel probability measure on $X$ satisfying \eqref{4-1}.
	  Assume in addition  that the  condition (b) in Theorem \ref{thm15} is satisfied. 
	 Then for any  $g\in L^\infty(X, \dd\mu)$ with
         $\|g\|_{L^\infty(\dd\mu)}\leq 1$,  and  each positive integer
         $N\ge 20$,   there exist points  $y_1, \ldots, y_{_{2(r + 2)
             N}} \in X$ and  real 
	numbers $\lambda_1, \ldots, \lambda_{_{2(r + 2) N}}$ such that
	\begin{align*}
	&  \max_{x \in X} \left| \int_X \Phi (\rho (x, y)) g(y) \,\dd \mu (y) - \sum_{j =
		1}^{2(r + 2) N} \lambda_j \Phi (\rho (x, y_j)) \right|   \leqslant 45c_3N^{- \frac{1}{2} - \frac{3}{2\be}}  \sqrt{\log N}.
	\end{align*}

\end{thm}

Let us give some  examples of the metric spaces $(X,\rho)$ and  the associated classes  $\mathcal{S}_{r}$ which satisfy the conditions of Theorem \ref{cor-4-2} .

\begin{exam}
	(i) Let $X=\sph$ be the unit sphere of $\RR^{d+1}$ equipped with the usual geodesic distance $\rho(x,y)=\arccos x\cdot y$ for $x, y\in \sph$.  If  $\vi\in C[-1,1]$ is a piecewise algebraic polynomial  of degree at most $n_0$ on $[-1,1]$, then  the function $\Phi(\ta):=\vi(\cos\ta)$, $\ta\in [0,\pi]$  belongs to a  class $\mathcal{S}_r$ with $r$ being the dimension of the space of all spherical polynomials of degree at most $n_0$ on the sphere $\sph$. In this case, $\Phi(\rho(x,y)) =\vi(x\cdot y)$, and  the condition \eqref{4-1} implies both the condition (a) and the condition (b).
	
	(ii) Let $X=B_{\f \pi2}(0)\subset \RR^d$ be the Euclidean  ball with centre $0$ and radius $\f \pi2$.  If  $\vi\in C[0,\infty)$ is a piecewise algebraic polynomial  of degree at most $n_0$ on $[0,\infty)$, then  the function $\Phi(t):=\vi(t^2)$, $t\ge 0$  belongs to a  class $\mathcal{S}_r$ with $r$ being the dimension of the space of all algebraic  polynomials of degree at most $2n_0$ in $d$ variables. In this case, $\Phi(\rho(x,y)) =\vi(\|x-y\|^2)$, and the condition \eqref{4-1} implies both the condition (a) and the condition (b).
\end{exam}

We will discuss these examples in detail in Sections 7 and 8.

\subsection{Proof of Theorem \ref{thm15}}
\noindent
The proof of Theorem \ref{thm15}
follows along the same idea as that of \cite{BL}.

Let $\{X_1, \ldots, X_N\}$ be a partition of $X$ satisfying the condition (a).  By the inner regularity of
the measure $\mu$, for each
$1 \leqslant j \leqslant N$, there  exists a compact subset $Q_j \subset X_j $  such that
\[ \frac{1}{N} - \mu (Q_j) \leqslant \frac12 (1 +\| \Phi \|_{\infty})^{- 1}
   N^{- \frac{3}{2} - \frac{3}{2\be } } . \]
Let $\mu_j$ denote the probability measure on $Q_j$ given by $\mu_j (E) =
\frac{\mu (E)}{\mu (Q_j)}$ for each Borel subset $E \subset Q_j$. Then it is
easily seen that
\begin{equation}
  \sup_{x \in X} \left| \int_X \Phi (\rho (x, y)) \,\dd \mu (y) - \frac{1}{N}
  \sum_{j = 1}^N \int_{Q_j} \Phi (\rho (x, y)) \,\dd \mu_j (y) \right|
  \leqslant N^{- \frac{1}{2} - \frac{3}{2\be } } .  \label{5-8}
\end{equation}
Let $\Sigma_j$ denote the set of all Borel probability measures $\sigma_j $ on
$Q_j$ that take  the form
\[ \sigma_j = \sum_{i = 1}^{r + 2} \lambda_i (\sigma_j) \delta_{y_i
   (\sigma_j)},\  \ \lambda_i (\sigma_j) \ge 0,\  \ y_i (\sigma_j) \in Q_j, \   \ 1
\leqslant j \leqslant r + 2,  \]
 such that $\sum_{i=1}^{r+2} \ld_i(\sa_j)=1$ and
\begin{equation}
  \int_{Q_j} f (y)  \hspace{0.17em} \,\dd \mu_j (y) =\sum_{i=1}^{r+2} \ld_i(\sa_j) f(y_i(\sa_j)), \  \  \forall f \in V_r . \label{5-10}
\end{equation}
According to Theorem \ref{thm-6-1}, there exists a Borel
probability measure $\nu_j$ on $\Sigma_j$ such that
\begin{equation}
   \int_{Q_j} f \hspace{0.17em} \,\dd \mu_j = \int_{\Sigma_j} \sum_{i
  = 1}^{r + 2} \lambda_i (\sigma_j) f (y_i (\sigma_j))  \hspace{0.17em} \dd
  \nu_j (\sigma_j),\  \quad \forall f \in C (Q_j) . \label{5-11}
\end{equation}
Now we  consider the following product  probability space:
\[ (\wt{\Sigma}, \nu) = \prod_{j = 1}^N (\Sigma_j, \nu_j) . \]

 We first claim that for  each fixed $x \in X$ and  parameter  $t > \sqrt{\log2}$, there exists a subset $G(x)\subset \wt{\Sigma} $ with  $\nu (G(x))\leq 2 e^{- t^2}<1$  such that for each  $\sigma:= (\sigma_1, \sigma_2, \ldots, \sigma_N) \in
\wt\Sigma\setminus G(x)$,
\begin{align}
 & \left| \frac{1}{N}  \sum_{j = 1}^N \sum_{i = 1}^{r + 2} \lambda_i (\sigma_j)
  \Phi (\rho ((x, \nobracket y_i (\sigma_j))) - \frac{1}{N} \sum_{j = 1}^N
  \int_{Q_j} \Phi (\rho (x, y)) \,\dd \mu_j (y) \right| \notag\\
 & \leqslant \f {4}{\sqrt{3}} c_1\sqrt{c_1c_2\ell} tN^{- \frac{1}{2} - \frac{3}{2\be}}.\label{4-7}
\end{align}
To show this claim, we consider the following independent random variables
on the probability space $(\wt{\Sigma} , \nu)$:
\[ h_j (\sigma) \equiv h_j (\sigma_j) \assign \sum_{i = 1}^{r + 2} \lambda_i
   (\sigma_j) \Phi \Bigl(\rho (x, y_i (\sigma_j))\Bigr) - \int_{Q_j} \Phi (\rho (x, y))
   \,\dd \mu_j (y), \]
   where  $\sigma = (\sigma_1, \ldots, \sigma_N)
   \in \wt\Sigma$  and $j = 1, \ldots, N.$
 By {\eqref{4-3-0}} and
{\eqref{5-11}}, we have
\[ \mathbb{E}h_j = 0, \quad |h_j | \leq \text{diam}(X_j)\leq  \delta_N,\   \ 1 \leq j \leq N. \]
For each $1\leq j\leq N$, pick a point $y_j\in Q_j$ and set  $R_j \assign B_{\delta_N} [y_j]$ so that   $Q_j \subset X_j \subset R_j$.
Set
$$ S_i(x):=E(x; t_{i-1}, t_i)=\Bl\{ y\in X:  t_{i-1} \leq \rho(x,y) \leq t_i\Br\},\  \quad i=1,\ldots, \ell.$$
Note that if  $R_j \subseteq S_k (x)$ for some $1 \leqslant k
\leqslant \ell$ and $1\leq j\leq N$, then there exists a function $f_{k,x}\in V_r$ such that 
$$\Phi (\rho (x, \cdot)) \Bl|_{Q_j}=f_{k,x}\Bl|_{Q_j},$$
 which,  using {\eqref{5-10}}, implies that
\[  h_j (\sigma_j) = \sum_{i = 1}^{r + 2} \lambda_i
(\sigma_j) f_{k,x}(y_i (\sigma_j)) - \int_{Q_j} f_{k,x} (y)
\,\dd \mu_j (y)=0.\]
For $1 \leqslant k \leqslant \ell - 1$ and if $\ell>1$, let
\[ E_k (x) \assign \left\{ y \in X : \; t_k - \delta_N \leqslant
   \rho (x, y) \leqslant t_k +  \delta_N \right\} . \]
Denote by $I$ the set of all positive integers $1 \leqslant j \leqslant N$
such that
\[ y_j\in   \bigcup_{k = 1}^{\ell - 1} E_k (x). \]
Let $I^c = \{ 1, 2, \ldots, N \} \setminus I$. Note that  if  $j \in I^c$, then there exists $1\leq k\leq \ell$ such that $R_j\subset S_k(x)$, which implies  $h_j = 0$. Furthermore,  since
\[ \bigcup_{j \in I} X_j \subseteq \bigcup_{j\in I} R_j\subseteq \bigcup_{k = 1}^{\ell - 1} \left\{ y \in X
   : \quad t_k -2\delta_N \leqslant \rho (x, y) \leqslant t_k + 2
   \delta_N \right\}, \]
it follows by Condition (b) that
\[ \#I \leq 2c_2\ell N\delta_N = 2c_2 c_1\ell N^{1 -
  \be^{-1}} . \]
We shall use this in our next estimate.
Now setting
\[ \xi_j = \frac{1}{\delta_N} h_j, \quad j = 1, 2, \ldots, N, \]
and using  the Bernstein inequality in probability, we obtain that for any $\varepsilon > 0$,
\begin{align*}
  \text{Prob} \left\{ \frac{1}{N} \left| \sum_{j = 1}^N \xi_j \right| >
  \varepsilon \right\} & = \text{Prob} \left\{ \frac{1}{\#I} \left| \sum_{j
  \in I} \xi_j \right| > \frac{\varepsilon N}{\#I} \right\}\\
  & \leq 2 \exp \left( - \frac{3}{8} (\#I) \frac{\varepsilon^2 N^2}{(\#I)^2}
  \right) \leq 2 \exp \left( - \frac{3\varepsilon^2 N^{1 + \be^{-1}} }{16 c_1c_2\ell}
  \right) .
\end{align*}
It follows that for any $\delta > 0$,
\begin{align*}
  \text{Prob} \left\{ \frac{1}{N} \left| \sum_{j = 1}^N h_j \right| > \delta
  \right\} \leq 2 \exp \left( - \frac{3\delta^2 N^{1 + 3 \be^{-1}}}{16 c_1^3
  c_2\ell} \right) .
\end{align*}
Given a parameter $t>0$, setting
\[ \delta :=\f 4{\sqrt{3}} c_1\sqrt{c_1c_2\ell}  N^{- \frac{1}{2} - \frac{3 }{2\be}} t, \]
we   conclude that the inequality
\[ \frac{1}{N} \left| \sum_{j = 1}^N h_j \right| \leq \f 4{\sqrt{3}}
c_1\sqrt{c_1c_2\ell}\cdot t\cdot N^{-
   \frac{1}{2} - \frac{3 }{2\be}} \]
holds with probability at least $1 - 2 e^{- t^2}$ on the probability space
$(\wt{\Sigma} , \nu)$. This proves the claim \eqref{4-7}.

  Now let  $t \assign \sqrt{A \log N}\ge \sqrt{\log 2}$ with  $A > 1$ being a parameter to be
  specified later. 
  By \eqref{5-8} and \eqref{4-7}, for each $x\in X$,
  there exists a set $G (x) \subset \wt{\Sigma} $ with $\nu (G
  (x)) \leq 2 N^{- A}$ such that for each
  $$\sigma = (\sigma_1, \ldots, \sigma_N) \in \wt{\Sigma}  \setminus G (x),$$
    \begin{align}
    \biggl|& \frac{1}{N}  \sum_{j = 1}^N \sum_{i = 1}^{r + 2} \lambda_i
    (\sigma_j) \Phi (\rho(x, y_i (\sigma_j))) - \Phi_0(x) \biggr|\notag\\
    & \leq \f 72 c_1\sqrt{c_1c_2\ell }
    \sqrt{A} \sqrt{\log N}  N^{- \frac{1}{2} - \frac{3}{2\be}},\label{4-8-0}
  \end{align}
  where
   \[\Phi_0(x):=\int_X \Phi (\rho (x, y)) \,\dd \mu (y).\]
  Let $M$ be a positive  integer such that \[ M-1< c_1^{\be} N^{\f 32+\f \be {2}} \leq M.\]
 Then, using Condition (a) with $M$ in place of $N$, we obtain a
 partition $$\{ X_1', X_2',\ldots, X_M'\}$$ of $X$ such that $\mu(X_j')
 =\f 1M$ and
 $$\text{diam}(X_j')\leq \da_M=c_1M^{-\be^{-1}} \leq N^{-\f12-\f3{2\be}}$$ for
 each $1\leq j\leq M.$
 Choose $z_j\in X_j'$ for each $1\leq j\leq M$, and let  $G = \bigcup_{k = 1}^M G (z_k)$. Then
 \[ \nu (G) \leq \sum_{j = 1}^M \nu (G (z_j)) \leq 2 MN^{- A} \leq 3c_1^{\be}
 N^{\frac{\be}{2} +\f32 - A}. \]
 Thus,  setting $A = \frac{1+2c_1}{2}\be +\f32$, we obtain
 that for  $N\ge 4$, $\nu (G)$ is at most
 $$
 3 c_1^{\be} N^{- {c_1\be}}\leq \biggl(\f{ 3 c_1 }{4^{c_1}}\biggr)^{\be} <1.$$

Finally, using \eqref{4-3-0}, we have that for each  $\sigma = (\sigma_1, \ldots,
\sigma_N) \in \wt{\Sigma}  \setminus G$
 \begin{align*}
    & \sup_{x \in X} \left| \frac{1}{N}  \sum_{j = 1}^N
    \sum_{i = 1}^{r + 2} \lambda_i (\sigma_j) \Phi (\rho(x,  y_i (\sigma_j)))
    - \Phi_0(x) \right|\\
    & \leq \max_{1 \leq k \leq M} \left| \frac{1}{N}  \sum_{j = 1}^N \sum_{i
    = 1}^{r + 2} \lambda_i (\sigma_j) \Phi (\rho(z_k,  y_i (\sigma_j))) -
    \Phi_0(x) \right| + \da_M,
  \end{align*}
  which, using \eqref{4-8-0},  is estimated from above by
  \begin{align*}
    & \biggl(\f 72 c_1^{\f 32} (c_2\ell)^{\f12} \sqrt{\f {2c_1+1}
      {2}\be +\f 32}+1\biggr)  N^{- \frac{1}{2} - \frac{3}{2\be}}
    \sqrt{\log N}\leq  8 c_1^{2} (c_2\ell)^{\f12}  {\sqrt{\be}}
    N^{- \frac{1}{2} - \frac{3}{2\be}}\sqrt{\log N}.
  \end{align*}
  This completes the proof.


  \subsection{Proof of Theorem \ref{cor-4-2}}
  \noindent
Let $h(x)=b(2+g(x))$, where $b$ is a normalizing constant so that $\|h\|_{L^1(\dd\mu)}=1$. Clearly, \begin{equation}\label{4-8}
\f 13\leq b\leq h(x) \leq 3b\leq 3,\  \ \forall x\in X,
\end{equation}
because $\|g\|_\infty\leq1$. Let $\tau$ denote the Borel  probability measure given by  $\dd\tau =h\dd\mu$.  By \eqref{4-1}, we have that for $N\ge 15$,
 \begin{equation}\label{4-9}
 \tau (B_{\wt{\delta}_N/8} (x)) \geq b\mu (B_{\wt{\delta}_N/8} (x))\ge   \frac{1}{N},\  \ x\in X,
\end{equation}
where $$\wt{\da}_N=c_1 ([Nb])^{-\be^{-1}}\leq \biggl( \f
5{4b}\biggr)^{1/\be} c_1 N^{-\be^{-1}}\leq \f 5{4b} c_1
N^{-\be^{-1}},$$
because $\beta\geq1$. Furthermore, by \eqref{4-1-0}, we have
that  for each $x \in X$ and $\da\in(0,\pi)$,
\begin{equation}
\tau \left( \bigcup_{j = 1}^{\ell - 1} \Bigl\{ y \in X : \ t_j -
\delta \leqslant \rho (x, y) \leqslant t_j +  \delta \Bigr\}
\right) \leqslant 3bc_2\ell \delta. \label{4-6}
\end{equation}
Since $X$ is a compact path-connected  metric space,  using  Theorem \ref{thm15} with  $\tau$ in place of $\mu$, we may find  points  $y_1, \ldots, y_{_{(r + 2) N}} \in X$ and nonnegative  real
numbers $a_1, \ldots, a_{_{(r + 2) N}},$ such that
\begin{align*}
&  \max_{x \in X} \left| \int_X \Phi (\rho (x, y)) h(y) \,\dd \mu (y) - \sum_{j =
	1}^{(r + 2) N} a_j \Phi (\rho (x, y_j)) \right|
   \leqslant \f{25\sqrt{3}} {16b^{\f 32}}c_3 N^{- \frac{1}{2} - \frac{3}{2\be}}  \sqrt{\log N}.
\end{align*}
On the other hand, using Theorem \ref{thm15}, we can also find points  $z_{1}, \ldots, z_{(r + 2) N} \in X$ and nonnegative  real
numbers $b_{1}, \ldots, b_{(r + 2) N},$ such that
\begin{align*}
&  \max_{x \in X} \left| \int_X \Phi (\rho (x, y)) \,\dd \mu (y) - \sum_{j =
	1}^{(r + 2) N} b_j \Phi (\rho (x, z_j)) \right|   \leqslant c_3  N^{- \frac{1}{2} - \frac{3}{2\be}}  \sqrt{\log N}.
\end{align*}
Since
\[ \int_{X}\Phi(\rho(x,y)) g(y)\,\dd \mu(y) =\f 1{b}\int_{X}\Phi(\rho(x,y)) h(y)\,\dd \mu(y)-2\int_{X}\Phi(\rho(x,y)) \,\dd \mu(y)\]
and $\f 13\leq b\leq 1$, it follows that \begin{align*}
\sup_{x\in X} &\Bl|  \int_{X}\Phi(\rho(x,y)) g(y)\,\dd \mu(y) - \f 1b  \sum_{j =
	1}^{(r + 2) N} a_j \Phi (\rho (x, y_j)) +2\sum_{j =
	1}^{(r + 2) N} b_j \Phi (\rho (x, z_j))\Br|\\
&\leq \biggl(  \f {25\sqrt{3}}{16b^{\f 52}}+ 2\biggr ) c_3 N^{- \frac{1}{2} - \frac{3}{2\be}}  \sqrt{\log N}\leq 45 c_3 N^{- \frac{1}{2} - \frac{3}{2\be}}  \sqrt{\log N}.
\end{align*}
The theorem is proved.

\section{Discretization  on finite-dimensional compact domains}
\noindent
In this section, we shall
 prove an analogue of  Theorem  \ref{cor-4-2}  for all  $g\in L^1(\dd\mu)$ (instead of $g\in L^\infty(\dd\mu)$)  on finite-dimensional domains. The implied  constant in this section will depend on the dimension and the  underlying domain.

Let $(X,\|\cdot\|)$ be a finite-dimensional real normed linear
space. Let $B_\zeta(x)$ (resp. $B_\zeta[x]$)  denote the open balls
(resp.\ closed balls) with centre $x\in X$ and radius $\zeta>0$
defined with respect to the metric $\rho(x,y)=\|x-y\|$. Here
$\|\cdot\|$ is not  necessarily the Euclidean norm. 
Let $\Og\subset B_1[0]$ be a   compact subset of $X$ (not necessarily connected). Let $\mu$  be    a Borel probability measure supported on $\Og$. 
 The main purpose in this section is to discretize integrals of the form 
$$\int_{\Og} \Phi(\|x-y\|) g(y)\,\dd \mu(y) \   \ \text{
   with}\   \ g\in L^1(\dd\mu)$$  for a class of  piecewisely defined
 functions $\Phi:[0,\infty)\to\RR$.

We assume that the probability measure $\mu$ satisfies     the following two conditions:

\begin{enumerate}[\rm (i)]
	\item  there exist a  positive  constant $c_4>1$ and a parameter $\beta\ge 1$  such that 
		 for any $x\in \Og$ and $\da\in (0,2]$
	\begin{equation}\label{6-1}
	c_4^{-1} \da^{\beta} \leq \mu\Bl(B_{\da} (x)\Br)\leq c_4 \da^{\beta};
	\end{equation}

	\item  there exists a constant  $c_5>0$ such that  for any  $x\in \Og$ and  $t, s\in (0,2]$, \begin{equation}\label{5-2-1}
	\mu\Bl(\{y\in\Og:\  t\leq \|y-x\|\leq t+s\}\Br)\leq c_5 s.
	\end{equation}

\end{enumerate}

Under these two conditions,  we shall prove 
\begin{thm}\label{thm-6-1} Let $\Phi:[0,\infty)\to\RR$ be a function such that 
	\begin{equation}
	\label{6-2-1} | \Phi (s) - \Phi (s') | \leq |s - s'|,\   \forall s, s' \in [0,
	2].
	\end{equation}
	Assume that  there exist a partition $0=t_0<t_1<\cdots<t_\ell=2$ of $[0,2]$ and a translation-invariant  linear subspace $X_r$ of $C(\Og)$  with $\dim X_r=r$ such that with  $E_j:=\{x\in\RR^d:\ \ t_{j-1}\leq \|x\|\leq t_j\}$,  $j=1,2,\ldots, \ell$,  \begin{align*}
	 \Phi(\|\cdot\|)\Bl|_{E_j}\in  \Bl\{f\Bl|_{E_j}:\  \ f\in
         X_r\Br\}.
        	\end{align*}
	 Let $g\in L^1(\Og, \mu)$ be such that $\|g\|_{L^1(\dd\mu)}=1$.
	Then for  each positive   integer $n\ge 2$,   there exist  points $y_1,\ldots, y_{n}\in \Og$ and
	real numbers $\ld_1,\ldots, \ld_{n},$ such that
	\begin{align}
	\sup_{x\in \Og}&	\Bl|\int_{\Og} \Phi(\|x-y\|) g(y)\,\dd \mu(y)-\sum_{k=1}^{n}  \ld_k \Phi(\|x-y_k\|) \Br|\notag\\
	& \leq C(X)  \begin{cases}
	n^{-\f12-\f 3{2\beta}} (\log n)^{\f 12}, &\  \ \text{if $1\leq \beta<3$},\\
	n^{-1} (\log n)^{\f 32}, &\  \ \text{if $\beta=3$},\\
			 n^{-\f{\beta+1}{2(\beta-1)}} (\log n)^{\f 12}, &\  \ \text{if $\beta>3$},	
	\end{cases}
	\end{align}
	where the constant $C(X)$ depends only on $\dim X$,
        $c_4$, $c_5$, $r$, $\ell$ and $\be$.
\end{thm}

\subsection{Proof of Theorem \ref{thm-6-1}}
The main idea of our proof comes from the paper \cite{BL}. We need    the following   Besicovitch covering theorem on finite-dimensional normed linear spaces \cite{Kr}:

\begin{lem} \cite{FL}  Let $E\subset X$ be an  arbitrarily given
	nonempty subset of a finite dimensional normed linear space $X$. Assume that for each $x\in E$ there exists a
	closed ball $B_{r(x)}[x]$ with centre $x$ and radius
	$r(x)>0$. Assume in addition that  $\sup_{x\in E} r(x)<\infty$.
	Then there exists a sub-collection $\mathcal{R}$ of the closed balls
	$B_{r(x)}[x]$, $x\in E$, which covers the set $E$ and can be written
	in the form
	$$\mathcal{R}=\mathcal{R}_1\cup\mathcal{R}_2\cup\cdots\cup\mathcal{R}_m$$
	with $m\leq \mathcal{N}(X)$, and each $\mathcal{R}_j$ being a
	collection of pairwise disjoint balls, $1\leq j\leq m$. Here
	$\mathcal{N}(X)$ is a positive constant  depending only on the normed
	space $(X,\|\cdot\|)$.
\end{lem}
The best constant $\mathcal{N}(X)$ for the Besicovitch covering
theorem has been well studied in literature (see \cite{FL, Kr} and
the references thererin). In the case when  $(X,\|\cdot\|)=\RR^d$, it
was  known \cite{Kr} that $\mathcal{N}(X)\leq 6^d$. The sharp estimate
of this constant   appears  in \cite{Su}. A much more general version
of the Besicovitch covering theorem can be found in  \cite{Fe}.

The proof runs along the same line as that of  Theorem \ref{thm15}. We sketch it as follows.

Without loss of generality, we may assume that $g\ge 0$ since
otherwise we may write $g=g^{+}-g^{-}$ with $g^\pm\geq0$.
  For the rest of the proof, the letter $C$  denotes a general
  positive constant  depending only on $\mathcal{N}(X)$, $c_4$, $c_5$, $r$, $\ell$ and $\beta$.

 Let  $\tau$ denote the probability measure  given by
$\dd\tau(x)=g(x)\,\dd \mu(x).$ Let $n_1=[\f n {2\mathcal{N}(X) (r+2)}]$.
For $x\in \Og$, let $0<\ta_x\leq \da_{n_1}:= (c_4/ n_1)^{\f 1\beta}$  be such that
\begin{equation}\label{6-3-0}
\int_{B_{\ta_x}[x]}(1+ g(y)) \,\dd \mu(y)=\f 1 {n_1}.
\end{equation}
By the Besicovitch covering theorem, we can find finitely many open balls $B_{j} =B_{\ta_{x_j}}(x_j)$, $j=1,2,\ldots, m,$
such that $\Og \subset \bigcup_{j=1}^m B_j$,
\begin{equation}\label{6-4-0}
 \{B_1,\ldots, B_m\} =\mathcal{R}_1\cup \mathcal{R}_2\cup\cdots\cup\mathcal{R}_{\mathcal{N}(X)}
\end{equation}
with each $\mathcal{R}_j$ being a subcollection of pairwise disjoint balls.
By \eqref{6-3-0} and \eqref{6-4-0}, we then have
$m\leq 2\mathcal{N}(X) n_1\leq \f n {r+2}$. Note that \eqref{5-2-1} implies that $\mu(B_r(x))=\mu(B_r[x])$ for any $x\in\Og$ and $r>0$.
Now   define
$Q_1=\overline{B_1}$ and
$$ Q_j=\overline{B_j} \setminus \bigcup_{i=1}^{j-1} B_i,\    \  j=2,\ldots, m.$$
Then $\Og=\bigcup_{j=1}^m Q_j$,   $\tau(Q_i\cap Q_j)=0$ for $1\leq i\neq j\leq m$, $Q_j\subset \overline{B_j}$ and  $\tau(Q_j) \leq \f 1 {n_1}$ for $1\leq j\leq m$. Without loss of generality, we may also assume that $\tau(Q_j)>0$ for each $1\leq j\leq m$, since otherwise we remove $Q_j$ from the partition.

For each $1\leq j\leq m$, let $\Sigma_j$ denote the set of all  probability measures $\sa_j$  on $Q_j$ of the form $$\sa_j=\sum_{i=1}^{r+2} \ld_i(\sa_j) \da_{y_i(\sa_j)},\   \ \ld_i(\sa_j)\ge 0,\   \  y_i(\sa_j)\in Q_j,$$
such that
$$ \f 1 {\tau(Q_j)}\int_{Q_j} P(x)\,\dd \tau(x) =\sum_{i=1}^{r+2} \ld_i(\s_j) P(y_i(\s_j)),\   \  \forall P\in X_r.$$
By Theorem \ref{thm-6-1},  there exists a Borel probability measure $\nu_j$ on $\Sigma_j$ such that
$$ \int_{Q_j} f(x)\,\dd \tau(x) =\int_{\Sigma_j} \sum_{i=1}^{r+2} \tau(Q_j)\ld_i(\sa_j) f(y_i(\sa_j))\,\dd \nu_j(\s_j),\   \  \forall f\in C(Q_j).$$
Now we consider the product probability space
$(\wt \Sigma, \nu)=\prod_{j=1}^m (\Sigma_j,\nu_j)$. Fix  $x\in \Og$  temporarily.
For $1\leq j\leq m$,   define
$$h_{j,x} (\sa_j) =\tau(Q_j) \sum_{i=1}^{r+2} \ld_i (\sa_j) \Phi(\|x-y_i(\sa_j)\|) -\int_{Q_j} \Phi(\|x-y\|) \,\dd \tau(y).$$
Then
$\EE h_{j, x} =0$, \begin{equation}\label{5-6-0}
|h_{j,x} (\sa_j)|\leq\tau(Q_j)\cdot \diam(Q_j)\leq \tau(B_j) \diam(Q_j)\leq  C\ta_{x_j} n^{-1}.
\end{equation}

\noindent
For $0<\ta\leq 2\da_{n_1}$, we denote by $I_\ta:=I_\ta(x)$ the set of all integers $1\leq j\leq m$ such that  $ \ta/2< \ta_{x_j}\leq\ta $ and  $t_k-\ta\leq \|x-x_j\|\leq t_k+\ta$ for some $1\leq k\leq \ell$.
Note  that   if  $ \ta/2< \ta_{x_j}\leq \ta$ and  $j\notin
I_{\ta}(x)$, then  there exists $k$ in the interval $1\leq k\leq \ell$
such that $t_{k-1}\leq \|x-y\|\leq t_k$ for every $y\in Q_j\subset
B_j:= B_{\ta_{x_j}} (x_j)$, which implies that $
h_{j,x}\equiv 0$. Note also that $$\bigcup_{j\in I_\ta} \Bl( B_j \cap \Og\Br) \subset \Bl\{ y\in\Og:\  t-2\ta\leq \|x-y\|\leq t+2\ta\Br\}.$$
It then follows by  \eqref{6-4-0} and \eqref{5-2-1} that
$$ \# I_\ta c_4 \biggl( \f \ta2\biggr)^{\beta} \leq \sum_{j\in I_\ta}\mu (B_j) \leq 4\mathcal{N}(X) c_5\ta,$$
which implies that \begin{align}\label{5-6}
\#I_\ta& \leq C_1\ta^{1-\beta}.
\end{align}
Note that \eqref{5-6} holds trivially
 if \begin{equation}
 \ta\leq \biggl( \f {(r+2)C_1}{n}\biggr)^{\f 1{\beta-1}} =C_2 n^{-\f 1{\beta-1}}
 \end{equation} since   $\# I_\ta\leq m\leq \f {n}{r+2}$.
Thus, we will mainly  consider those index sets $I_\ta$ with
\begin{equation}
 C_2	n^{-\f 1{\beta-1}} \leq  \ta\leq 2 \da_{n_1}:= 2\Bl(\f
 {c_4}n\Br)^{\f 1\beta}, 
\end{equation}
the second bound being the bound on $\theta$
stated at the beginning of the paragraph.

To be more precise, let  $k_0$, $k_1$  be  integers such that
$$ 2^{k_0}<2^{-1}(n/c_4)^{\f1\beta}\leq 2^{k_0+1}$$
and $$2^{k_1-1}<  C_2^{-1} n^{\f 1{\beta-1}}\leq 2^{k_1}.$$	
Define $J_k=J_k(x):=I_{2^{-k}}(x)$ for $k_0\leq k\leq k_1$ and
$$J_{k_1+1}\equiv J_{k_1+1} (x)=\bigcup_{k=k_1+1}^\infty
I_{2^{-k}}(x).$$
Then by \eqref{5-6} and the remark after \eqref{5-6}, we have
\begin{equation}\label{5-9}
\#J_k \leq n_k:= C_1^{-1} 2^{k (\beta-1)},\   \ k_0\leq k\leq k_1+1.
\end{equation}
Moreover, by \eqref{5-6-0}, we have \begin{equation}\label{5-11-0}
|h_{j,x}|\leq C \ta_{x_j} n^{-1}\leq C 2^{-k} n^{-1},\   \ j\in J_k,\  \ k_0\leq k\leq k_1+1.
\end{equation}
Thus,  using \eqref{5-11-0}, \eqref{5-9}, and the Bernstein inequality, we  conclude that for  each   $k_0\leq k\leq k_1+1$ and each $\va_k>0$,  the inequality
\begin{align*}
\Bl|\sum_{j\in J_k}  h_{j,x}(\sa_j)\Br| >\va_k
\end{align*}
holds with probability at most
\begin{equation}\label{5-8-1}
2 \exp \Bl( -C \va_k^2 n^{2}2^{-k(\beta-3)}\Br).
\end{equation}

\noindent
Now we write
\begin{align*}
\sum_{j=1}^m h_{j,x} (\sa_j) &=\sum_{k=k_0}^\infty \sum_{\{j: 2^{-k}\leq \ta_{x_j}\leq 2^{-k+1}\}} h_{j,x}(\sa_j)=\sum_{k=k_0}^{k_1+1} \sum_{j\in J_k} h_{j,x}(\sa_j).
\end{align*}
  Given $\va>0$, let $\{\va_k\}_{k=k_0}^{k_1+1}$ be a sequence of positive numbers  such that $\sum_{k=k_0}^{k_1+1} \va_k \leq \va$. Then using \eqref{5-8-1}, we have
\begin{align}
\pr \Bl\{ |\sum_{j=1}^m h_{j,x}| >\va\Br \}&\leq \sum_{k=k_0}^{k_1+1}\pr \Bl\{ |\sum_{j\in J_k}  h_{j,x}| >\va_k \Br \}\notag\\
&\leq 2 \sum_{k=k_0}^{k_1+1}\exp \Bl ( -C \va_k^2 n^{ 2} 2^{-k(\beta-3)}\Br).\label{5-14}
\end{align}
Noting  that $k_0\sim k_1\sim \log n$, we may choose for $k_0\leq k\leq k_1+1$,
$$\va_k=\begin{cases} 2^{\f {\beta-3}2(k-k_1)} \va,&\  \ \text{if $\beta>3$},\\
 \f \va{\log n}, & \  \ \text{if $\beta=3$},\\
 2^{(k-k_0)\f {\beta-3}2}\va, & \  \ \text{if $\beta <3$}.
\end{cases} $$
\noindent
We use here that $n\neq1$ so that $\log n\neq0$.

For simplicity, we shall assume that $\beta>3$. The proof below with slight modifications works equally well for the case $\beta\leq 3$.
We then  obtain from \eqref{5-14} that
\begin{align*}
\pr \Bl\{ \Big|\sum_{j=1}^m h_{j,x}\Big| >\va\Br \}&\leq  C (\log n) \exp \Bl( - C n^{2-\f{\beta-3}{\beta-1}}\va^2\Br).
\end{align*}
Setting
$$\va=C^{-\f12} t\cdot n^{-\f{\beta+1}{2(\beta-1)}}\   \  \ \text{with}\  \ t>0,$$
we conclude that for each $x\in \Og$,
the inequality
\begin{align*}
&\Bl|\sum_{j=1}^m  \tau(Q_j) \sum_{i=1}^{r+2} \ld_i (\sa_j) \Phi(\|x-y_i(\sa_j)\|) -\int_{\Og} \Phi(\|x-y\|) \,\dd \tau(y)\Br|\\
&\ge
C t n^{-\f{\beta+1}{2(\beta-1)}}
\end{align*}
holds with probability bounded above by a multiple of
$ (\log n) e^{- t^2}$.
Let
$$ t:=\sqrt{A\log n}\  \ \text{with }
\  \ A=\frac{\beta(\beta+1)}{2(\beta-1)}>1.$$
The last inequality holds since $\beta^2-2\beta+2$ has no real
zeos.

We further conclude that for each $x\in\Og$, there exists a set $G(x)\subset \Sigma$ with $\nu(G(x)) \leq C_2(\log n) n^{-A} $ such that for any $\sa=(\sa_1,\ldots,\sa_m)\in\Sigma\setminus G(x)$,
\begin{align*}
&\Bl|\sum_{j=1}^m  \tau(Q_j) \sum_{i=1}^{r+2} \ld_i (\sa_j) \Phi(\|x-y_i(\sa_j)\|) -\int_{\Og} \Phi(\|x-y\|) \,\dd \tau(y)\Br|\\
&\leq
C^{-1/2} \sqrt{A}  n^{-\f{\beta+1}{2(\beta-1)}} (\log n)^{\f 12}.
\end{align*}

Finally, let $\{z_1, \ldots, z_L\}$ be a maximal $\va_1$-separated
subset of $\Og$ with $\va_1:=n^{-\f{\beta+1}{2(\beta-1)}} (\log n)^{\f
  12}$. By \eqref{6-1}, we have  \textbf{$$L\leq c_4\biggl( \f 2
{\va_1}\biggr)^\beta\leq C_3 n^{\f{\beta(\beta+1)}{2(\beta-1)}} (\log n)^{-\f
  12\beta}.$$}
Setting
$ A=\f{\beta(\beta+1)}{2(\beta-1)},$
we have that\textbf{\begin{align*}
\sum_{j=1}^L \nu (G(z_j)) \leq C_2 C_3   (\log n)^{1-\f 12\beta}.
\end{align*}}
Since $\be>3$, it follows that the following inequality holds with positive probability:
\begin{align*}
&\sup_{x\in\Og} \Bl|\sum_{j=1}^m  \tau(Q_j) \sum_{i=1}^{r+2} \ld_i (\sa_j) \Phi(\|x-y_i(\sa_j)\|) -\int_{\Og} \Phi(\|x-y\|) \,\dd \tau(y)\Br|\\
&\leq
C  n^{-\f{\beta+1}{2(\beta-1)}} (\log n)^{\f 12}.
\end{align*}
The theorem is proved.

\section{Discretization on  the unit sphere $\SS^d$}

\noindent
In this section, we will estimate the constants $c_1$ and $c_2$ for the
unit sphere  $\mathbb{S}^d\subset \RR^{d+1}$ denote the unit sphere in $\mathbb{R}^{d + 1}$ equipped
with the normalized  surface Lebesgue measure $\mu_d$  and the geodesic distance  $\rho (x, y) = \arccos (x \cdot y), x, y \in \mathbb{S}^d
.$ We will prove on the unit sphere $\SS^d$ that \begin{equation}\label{5-1-0}
c_1\leq 40 \pi,\  \  c_2\leq \f 32 \sqrt{d},\  \  \al=\f1d.
\end{equation}
The main point here lies in the fact that the upper bounds for $c_1$ and $c_2/\sqrt{d}$  are   independent of the dimension $d$.

By \eqref{5-1-0}, we also have \begin{equation}
45c_3=45\cdot 8c_1^2 c_2^{\f12} \sqrt{d}\leq 7\times 10^6 d^{\f 34}.
\end{equation}

As a consequence of Theorem \ref{thm10} and Lemma \ref{lem-5-2}, we have that
\begin{thm}
	\label{thm-1-3} For each integer $N\ge 1$, there exists a partition $\{R_1,
	\ldots, R_N \}$ of $\mathbb{S}^d$ such that
	\begin{enumerate}[\rm (i)]
		\item the $R_j$ are pairwise disjoint subsets of $\mathbb{S}^d$;
		
		\item for each $1 \leq j \leq N$, $\mu_d (R_j) = \frac{1}{N}$ and
		\text{diam}$(R_j) \leq 40
		\pi N^{- \frac{1}{d}}$.
	\end{enumerate}
\end{thm}

Again, the main point here is that the upper bound for $N^{\f 1d} \max_j \text{diam}(R_j)$ is independent of the dimension $d$.

\begin{thm}\label{thm-5-2}
	Let  $\Phi:[-1,1]\to\RR$   be  a piecewise polynomial of degree at most $r$   with  knots $-1=s_0<s_1<\cdots<s_\ell=1$ such that  $|\Phi(s)-\Phi(s')|\leq |s-s'|$ for any $s,s'\in [-1,1]$. Let $m_r=m_r^d$ denote the dimension of the space of all spherical polynomials of degree at most $r$ on $\SS^d$.
	 Let $g\in L^\infty(\sph)$ be such that $\|g\|_\infty \leq
         1$. Then for each positive integer $N\ge 20$, there exist points $\xi_1,\ldots, \xi_{2(m_r+2)N}\in\sph$ and real  numbers $\ld_1,\ldots,\ld_{2(m_r+2)N}$ such that
	 \begin{align*}
	 	 \max_{x\in\sph}& \left| \int_{\mathbb{S}^{d }} \Phi(x\cdot y) g(y)\,\dd \mu_d  (y) - \sum_{j = 1}^{2(m_r+2)N}
	 \lambda_{j} \Phi(x\cdot \xi_j) \right|\\
	 & \leq 7\cdot 10^6 \sqrt{\ell} d^{\f34} N^{- \frac{1}{2} - \frac{3}{2 d}}\sqrt{\log N} .
	 \end{align*}
\end{thm}
In the case when $\Phi(t)=|t|$, Theorem \ref{thm-5-2}, but with  constants depending on the dimension of the sphere, was previously obtained in \cite{BL}.

\subsection{Proof of \eqref{5-1-0}}

\noindent
For $\theta \in (0, \pi)$ and $x$ in the $d$-dimensional sphere
$\mathbb{S}^d$, set
\[ B_\ta (x) \assign \{y \in \mathbb{S}^d : \rho (x, y) <
   \theta\},\  \ \text{and}\  \  B_\ta [x] \assign \{y \in \mathbb{S}^d : \rho
   (x, y) \leq  \theta\} . \]
Let
$ \omega_d := \frac{2 \pi^{\frac{d + 1}{2}}}{\Gamma (\frac{d + 1}{2})} $
denote the surface area of $\mathbb{S}^d$. Using  the following known
estimates on gamma functions \cite{AS},
\[ x^{1 - s} < \frac{\Gamma (x+1)}{\Gamma (x + s)} < (x + 1)^{1 - s},\quad x >
0,\quad s \in (0, 1), \]
we  have that
\begin{equation}\label{5-1}
 \pi^{- \frac{1}{2}} \left( \frac{d - 1}{2} \right)^{\frac{1}{2}} \leq
\frac{\omega_{d - 1}}{\omega_d} = \frac{\Gamma (\frac{d + 1}{2})}{\Gamma
	(\frac{d}{2}) \sqrt{\pi}} \leq \pi^{- \frac{1}{2}} \left( \frac{d + 1}{2}
\right)^{\frac{1}{2}} .
\end{equation}

\begin{lem}
	\label{lem-5-1}For $0 < \theta \leq \frac{\pi}{4}$ and $x \in \mathbb{S}^d$,
	\[ \frac{1}{
          \sqrt{2d}} \leq \frac{\mu_d (B_\ta(x))}{\sin^d \theta}
	\leq \frac{2}{\sqrt{d}} . \]
\end{lem}

\begin{proof}
	For $\ta\in (0, \pi]$, we have
		\begin{align*}
	\mu_d (B_\ta(x)) & = \frac{\omega_{d - 1}}{\omega_d}  \int_{\cos
		\theta}^1 (1 - t^2)^{\frac{d - 2}{2}}  \hspace{0.17em} \dd t =
	\frac{\omega_{d - 1}}{\omega_d}  \int_0^{\sin^2 \theta} t^{\frac{d -
			2}{2}}  (1 - t)^{- \frac{1}{2}}  \hspace{0.17em} \dd t.
	\end{align*}
	If $0 < \theta \leq \frac{\pi}{4}$, then for any $0 \leq t \leq \sin^2
	\theta$, we have
	\[ 1 \leq (1 - t)^{- \frac{1}{2}} \leq \sqrt{2} . \]
	Thus,
		\begin{align*}
	\frac{\omega_{d - 1}}{\omega_d}  \frac{2}{d} (\sin \theta)^d \leq \mu_d (A
	(x, \theta)) \leq \frac{\omega_{d - 1}}{\omega_d}  \frac{2 \sqrt{2}}{d}
	(\sin \theta)^d,
	\end{align*}
	which, using \eqref{5-1}, implies that 	
	\begin{align*}
	  \frac{1}{
            2\sqrt{d}} \leq \sqrt{\frac{1}{\pi}} d^{- \frac{1}{2}} \leq \sqrt{2}
	\pi^{- \frac{1}{2}} d^{- 1}  (d - 1)^{\frac{1}{2}} \leq \frac{\mu_d (A (x,
		\theta))}{\sin^d \theta} \leq 2 \pi^{- \frac{1}{2}} d^{- 1}  (d +
	1)^{\frac{1}{2}} \leq \frac{2}{\sqrt{d}}.
	\end{align*}

\end{proof}

The following lemma shows that $c_1\leq 40 \pi$:
\begin{lem}\label{lem-5-2}For any positive integer $N$, \begin{equation}\label{5-2}
\inf_{x\in\sph} 	\mu_d (B_{\da_N}(x))\ge \frac 1N\   \ \text{with}\  \ \da_N:=5\pi N^{-\f 1d}.
	\end{equation}
	
\end{lem}

	\begin{proof} We consider the following two cases:\\
		
	  {\it Case 1.}    $N\ge 2^{\f d2+1}\sqrt{d}$.
          \\

	 In this case, set
\[ \delta \assign \min \left\{ \theta: \quad 0\leq \ta\leq \f \pi4,\  \  \frac{1}{2 \sqrt{d}}
\sin^d  {\theta} \geqslant \frac{1}{N} \right\}  \]
and our condition on the $N$ ensures that $\delta$ be well-defined.
Using \ Lemma \ref{lem-5-1},  we have  that
\[ \mu_d (B_\da(x)) \geqslant \frac{1}{N}, \quad \forall x
\in \mathbb{S}^d . \]
It remains to estimate the constant $\delta .$ By definition of $\delta$, we have that
\[ \frac{1}{2 \sqrt{d}} \sin^d  {\delta} \geqslant \frac{1}{N} >
\frac{1}{2 \sqrt{d}} \sin^d  \frac{\delta}{2} . \]
This implies that
\[ \delta \leqslant  \pi \sin \frac{\delta}{2} <
\pi \left(
\frac{2 \sqrt{d}}{N} \right)^{\frac{1}{d}} \leqslant 2 \pi \ee^{\frac{1}{2
    \ee}} N^{- \frac{1}{d}}<3\pi N^{- \frac{1}{d}} . \]

Here we have used the fact that the maximum of $(\log y)/y$ is attained at $y=\ee$.	

	{\it Case 2.}    $1\leq N< 2^{\f d2+1}\sqrt{d}$.\\

	In this case,
	$$N^{-\f1d} > 2^{-\f 1d -\f 12} d^{-\f 1{2d}}\ge 2^{-\f 32}\ee^{-\f 1 {2\ee}}>0.2,$$
	and
	$$ \da_N =5\pi N^{-\f 1d} \ge\pi. $$

	Hence,  \eqref{5-2} holds trivially in this case.
	
	\end{proof}

The following lemma shows that $  c_2\leq \f 32 \sqrt{d}$:
\begin{lem} For any $\da>0$,  $x\in\sph$ and $t\in (0,\pi)$,
\begin{equation}
\mu_d \left(  \Bigl\{ y \in \sph : \ t -
\delta \leqslant \rho (x, y) \leqslant t +  \delta \Bigr\}
\right) \leqslant \f 32 \sqrt{d} \delta.
\end{equation}
\end{lem}
\begin{proof}Without loss of generality, we may assume that $0<t\leq \f \pi2$. Setting \[ S_\da(x):=\Bigl\{ y \in \sph : \ t -
	\delta \leqslant \rho (x, y) \leqslant t +  \delta \Bigr\}, \]
 and using \eqref{5-1},   we have\begin{align*}
	\mu_d (S_\da(x))
		&=\f{\o_{d-1}}{\o_d}\int_{\max\{t-\da,0\}}^{t+\da}\sin^{d-1} u\,\dd u
		\leq \pi^{-\f12} \Bl( \f {d+1}2\Br)^{\f12}2\da\\
		&\leq \f {2}{\sqrt{\pi}}\sqrt{d}<\f 32 \sqrt{d}.	\end{align*}
	\end{proof}
\section{Further Examples}
\noindent
Further examples for our results stem from the fact that not only
piecewise polynomials are suitable for our spaces $V_r$ of dimension $r$, but also
piecewise exponentials \cite{AR}, \cite{CCR} and \cite{CGR}, as well
as radial basis functions of compact support \cite{B1}, \cite{B2} and
\cite{JM1}, \cite{JM2}.

All these function spaces are defined not over piecewise polynomials
(splines) with a simple continuity condition, but for  instance over
piecewise exponentials.

In the most general form, see \cite{AR}, the exponential splines of
compact support are, say, in $d$ dimensions of degree $n-1$ for {\it
  equally spaced knots\/} defined as
distributions $B$ that satisfy
$$ B(\varphi)=\int_{[0,1]^n}\varphi(\Xi t)\exp(\lambda\cdot t)\,\dd t,$$
where $\varphi$ is a test-function from the Schwartz space $S$,
$\Xi$ is a linear map $\R^n\to\R^d$ and $\lambda $ is a vector from
$\R^n$ to define the exponentials. Alternatively we can write for
$\phi\in L^1_{loc}(\R^d)$
$$\int_{\R^d} B(x)\phi(x)\,\dd x=
\int_{[0,1]^n}\phi(\Xi t)\exp(\lambda\cdot t)\,\dd t.$$

In the multivariate setting, these functions are called exponential
box-splines, in the univariate case they are exponential
B-splines. The piecewise polynomial case
corresponds to $\lambda\equiv0$. They may also be conveniently defined
by their Fourier transforms
$$\prod_{j=1}^n\frac{\exp(\lambda_j-\i \xi_j\cdot x)}{\lambda_j-\i \xi_j\cdot x}.$$
Here, $\lambda=(\lambda_j)_{j=1}^n$ and $\Xi=(\xi_j)_{j=1}^n$.

The univariate piecewise polynomial case corresponds to
$\Xi=(1,1,\ldots,1) \in\R^n$, $d=1$.
In this case the splines are defined over the interval or cube for
$d=1$ and $d>1$, respectively, $\Xi[0,h]^n$, e.g., $h=1/\ell$ as in our
cases. The $V_r$ space is here the space of univariate exponential
splines spanned by the exponential B-splines with $d+1$ knots.

More generally, we can define {\it space of piecewise exponentials\/}
including piecewise polynomials and exponentials as the span of
$$x^{r_1}\exp(\lambda_ix), \quad r_i=0,1,\ldots,\tau_i-1,\quad i=1,2,\ldots n,$$
on each subinterval between two knots, now of no longer necessarily equally spaced
knots, of
dimension $r=\sum_{i=1}^n\tau_i$ when they are required to be continuous. The
$\lambda_i$s may be complex and must be pairwise distinct.

Special cases \cite{CCR} are $\lambda\equiv0$ (piecewise polynomials),
$\lambda\in i\R$ ($V_r$ containing piecewise trigonometric functions $\sin$, $\cos$ and
constants)
and $\lambda\in\R$, $V_r$ containing $\sinh$ and $\cosh$ and constants. In fact, it is
usual (but not necessary) to restrict the exponents that form the components of $\lambda$
to $\R\cup\i\R$. Examples for the spaces are the polynomials for some fixed maximal
degree (classical spline case) or the spans of, e.g.,
$$1,\cos(\Im\lambda t),\sin(\Im\lambda t),t\cos(\Im\lambda t),t\sin(\Im\lambda t),$$
or
$$1,\cosh(\Re\lambda t),\sinh(\Re\lambda t),t\cosh(\Re\lambda
t),t\sinh(\Re\lambda t).$$

These two examples are the suitable generalisations of the
$\Phi(t)=|t|$ case (piecewise linears)
referred to in the paragraph after the statement of Theorem
\ref{thm-1-3}. For higher powers, larger $r$ and more exponentials,
the other piecewise polynomials used in the first sentence of the
statement of  Theorem \ref{thm-1-3} are generalised.

Univariate piecewise polynomial B-splines on equally spaced knots can be generate in a computational useful, recursive way by
convolutions \cite{Ca} but now, for exponential splines we get a weight function, so that, for
the B-spline of degree $n$, the exponential spline
$$ \ee^{\lambda_i t}H\Bigl( t\Bigr)-\ee^{\lambda_i}H\Bigl(t-1\Bigr)\ee^{\lambda_i(t-1)}$$
needs to be convolved with itself $n$-times, once for the case of
piecewise linears multiplied with exponentials. In the display, $H$ denotes the
Heaviside function which is identically zero for negative argument and identically one for positive argument.

This results from the identities which we stated already in $s$ dimensions
$$ B\star f=\int_{[0,1]^n}\exp(\lambda\cdot t)f(\cdot-\Xi t)\,\dd t$$
or
$$B=\int_0^1\exp(\lambda_\gamma\underline t)\tilde B(\cdot-\xi_\gamma \underline t)\,\dd \underline t.$$
Here $B$ is the exponential box-spline as above, $\tilde B$ is the same with the direction
$\xi_\gamma$ removed from $\Xi$.

As with the piecewise polynomials and the special case of piecewise constants above, we consider the special case of piecewise exponentials only (no polynomials as in our example with $\sin,\cos,\cosh,\sinh$).

For this, consider again the vector of exponents $\lambda$, set $n=d$
and let $\tilde \lambda=\lambda\Xi^{-1}$. Then the spline is
$$ B(x)=\frac1{|\det\Xi|}\exp(\tilde\lambda\cdot x)\chi_{(0,1]^d}\bigl(\Xi^{-1}x\bigr),\qquad x\in\R^d.$$
  Here, $\chi$ is the characteristic function. Starting from this
  piecewise ``constant'' function (i.e., one that contains no
  polynomials, just one exponential),
  {\it other\/}  splines can
  be generated recursively by
  $$B(x)=\ee^{\mu\cdot x}\int_0^1\tilde B(x-\underline t
  \xi)\,\dd \underline t,$$
  where $B$ is the exponential spline with one direction $\xi$ more in
  the direction set and the $\mu$s are chosen arbitrarily from $\R^n$.

The corresponding {\it radial basis functions\/} of compact support
with exponentials are
$$\Bigl(1/\ee-\exp(-x)\Bigr)_+^\nu$$
and $$\Bigl(1-\exp(-(1-x)^\nu_+)\Bigr)^\mu$$
which are positive definite for suitable parameters $\mu$ and $\nu$
depending on the dimension because they are logarithmically monotone
of  order $\mu$ in the first case and of order $\min(\mu,\nu)$ in the
second case \cite{JM1}. The $V_r$s are then defined by the translates
$$(1/\ee-\exp(-|x|))_+^\nu$$
and $$(1-\exp(-(1-|x|)^\nu_+))^\mu,$$
respectively.


\end{document}